\documentclass[a4paper,12pt]{article}
 \usepackage[cp1251]{inputenc}
\usepackage[english]{babel}
\usepackage[tbtags]{amsmath}
\usepackage{amsfonts,amssymb,mathrsfs,amscd}
\usepackage{color}
\usepackage{graphics}
\usepackage{wrapfig}
\usepackage{euscript}
\usepackage[dvips]{graphicx}

\date{}
\newtheorem{theorem}{Theorem}[section]
\newtheorem{lemma}{Lemma}[section]

\newtheorem{definition}{Definition}[section]

\begin{document}
\renewcommand{\thesection}{\arabic{section}}
\renewcommand{\theequation}{\thesection.\arabic{equation}}
\csname @addtoreset\endcsname{equation}{section}
\large

\begin{center}\bf On approximation by tight wavelet frames on Vilenkin groups
\end{center}
\centerline{\bf Sergey Lukomskii\footnote{Sergey Lukomskii:\ LukomskiiSF@info.sgu.ru}, Iuliia Kruss\footnote{Iuliia Kruss:\ krussus@gmail.com}, Alexandr Vodolazov\footnote{Alexandr Vodolazov:\ vam21@yandex.ru}}
Affiliation: Department of Mathematics and Mechanics, Saratov State University, Saratov, Russian Federation.

\begin{abstract}
  We  consider  the approximate properties of tight  wavelet frames  on Vilenkin  group $G$.   Let $\{G_n\}_{n\in \mathbb{Z} }$ be a main chain of subgroups, $X$ be a set of characters.  We define a step function $\lambda({\chi})$ that is constant on cosets  ${G}_n^\bot\setminus{G}_{n-1}^\bot$ by equalities
$\lambda ({G}_n^\bot\setminus{G}_{n-1}^\bot)=\lambda_n>0$
for which  $\sum\frac{1}{\lambda_n}<\infty$.
 We find the order of approximation of functions $f$ for which
 $\int_X|\lambda( {\chi})\hat{f}(\chi)|^2d\nu(\chi)<\infty$. As a corollary, we obtain an approximation error  for functions from Sobolev spaces with logarithmic weight.

     Bibliography: 18 titles.
\end{abstract}
\noindent
 keywords: tight  wavelet frame, zero-dimensional groups, local fields, refinable function, trees.
\section*{Introduction}
 Tight wavelet frames (often called Parseval frames) on a straight line are an important tool in image processing \cite{HZ}. Tight wavelet frames based on MRA can be considered as a generalization of the MRA-based  orthonormal  wavelets. Publication \cite{RSh97}
in 1997, the principle of unitary extension for constructing tight wavelet frames in $\mathbb R$ initiated to a new wave of both theoretical research and applications to information processing. The frames obtained in this way lead to fast algorithms for both decomposition and reconstruction, which is the reason for their attractiveness. An example is tight frames constructed from centered $B$-splines of order $m$.
 I.Daubechies, B.Han, A.Ron, Z.Shen \cite{DHRS}  considered the problem of approximation by tight wavelet frames $(\psi^{(j)}_{k,h})$ and proved that if a function $f$ belongs to the Sobolev space $W_2^{m_1}(\mathbb R)$  with a power weight function  $(1+|\xi|)^{2m_1}$, then
$$
\|f-\sum_{k<n}\sum_{j=1}^r\sum_{h\in \mathbb Z} \langle f,\psi^{(j)}_{k,h}\rangle \psi^{(j)}_{k,h}\|_{L_2(\mathbb R)}=O(2^{-nm_1}).
$$
 A large list of publications on this subject can be found in \cite{HZ} and \cite{ChO}

Yu.Farkov \cite{YuF3} constructed  the first examples of tight frames on the Cantor dyadic group. Yu.Farkov, E.Lebedeva, and M.Skopina \cite{FLS}  gave a complete de\-scription of all tight  wavelet frames generated by Vilenkin polynomials and proposed an algorithm for finding the corresponding wavelet functions. In the same paper, the approxi\-mation properties of tight  frames are considered. An approximation error for functions from the Sobolev space with   power weight function is obtained.
 In \cite{YuF4} three types of compactly supported wavelet frames associated with Vilenkin functions  are described: 1) MRA-based tight frames, 2) frames, obtained from the Daubechies-type “admissible condition”, and 3) frames, based on the Vilenkin Parseval type kernels.  In \cite{YuF6} finite tight frames for the space $l^2({\mathbb Z})$ are described.

 In this  article we devote  the construction of tight  wavelet frames
  on Vilenkin groups without using the unitary  extension  principle and  study  the approximation properties of the constructed frames. We obtain an approxi\-mation error for functions from the Sobolev space with    logarithmic weight.

 The article is organized as follows.
In Sec. 1, we give the necessary information about locally compact zero-dimensional groups  and prove some auxiliary results.  The main role in our constructions and proofs  is played by Lemmas 1.4 and 1.5, which are valid on any zero-dimensional group.

 In Sec. 2, we discuss the problem of constructing a step refinable  function $\varphi\in \mathfrak{D}_{\mathfrak{G}_M}({\mathfrak{G}_{-N}})$ (\cite{LS2}). We solve this problem in the Fourier domain  and construct the Fourier transform of the refinable  function as the product   of the masks dilations. To find the mask, we use a tree whose nodes contain the values of the mask on cosets. Using these trees, we specify a way to construct all  masks that generate the refinable  function $\varphi\in \mathfrak{D}_{\mathfrak{G}_M}({\mathfrak{G}_{-N}})$.

  In section 3, we indicate a method for constructing step tight wavelet frames on zero-dimensional group. This method is inde\-pendent of a particular group.   If the mask is constructed, then the procedure for constructing tight frames is the same for all zero-dimensional groups.

  In  section 4  we discuss  the approximation properties of tight  wavelet frames. We prove a general theorem 4.1, from which we obtain specific assertions.

\section{The Vilenkin  group as a locally compact zero-dimensional abelian group}
Let $p\ge 2$ be a fixed prime number. The Vilenkin group $G$ is the set all sequences
$$
x=(...,0,0,x_n,x_{n+1},x_{n+2},...), n\in \mathbb Z
$$
where $x_k\in \{0,1,...,p-1\}$ and there exists at most a finite number $k$ for  which $x_k\neq 0$. The group operation $\dot+$ is defined as the coordinatewise addition modulo $p$. Subgroups
$$
G_n=\{x=(...0,0_{n-1},x_n,x_{n+1},...  ):n\in\mathbb Z\}
$$
form the base of the topology. It is evident  $G_{n+1}\subset G_n,\  \bigcup_{n=-\infty}^{+\infty}G_n= G,$ $
 \quad \bigcap_{n=-\infty}^{+\infty}G_n=\{0\},$
  and $$G_n=\bigsqcup_{\nu=0}^{p-1}\bigl(G_{n+1}\dot+\nu g_n\bigr)
$$
where $g_n=(...,0,0_{n-1},1_n,0_{n+1})$. It is clear that $pg_n=0=(0)$. It means that
 $(G,\dot + )$~is a~locally compact zero-dimensional abelian
group with the topology generated by a~countable system of open
subgroups $G_n$. \cite{AVDR}

Therefore, we will carry out all reasoning in a locally compact zero-dimensional group  $G$ with the condition $pg_n=0$ on the operation $\dot+$. In this case $g_n\in G_n\setminus G_{n+1}$ and he system $(g_n)_{n\in \mathbb Z}$ is  called the basic system, since
  any element $g\in G$ has a~unique representation
in the form
\begin{equation}
\label{eq1.01}
g=\sum_{n=-\infty}^{+\infty}a_ng_n, \qquad
a_n=\overline{0,p-1}.
\end{equation}
The sum (1.1) contains a finite number of terms with negative numbers. Let $\mu $ be a Haar measure on $G$ and $\int_E fd\mu$ be an  absolutely convergent integral over a measure $\mu$.

 The mapping $\lambda:G\to[0,+\infty)$ defined by the equality
$$
\lambda(g)=\sum_{n=m}^{+\infty}a_np^{-n-1}
$$
is called Monna mapping \cite{Mon}.
 It is evident that $\lambda(G_n)=[0,p^{-n}]$.
 Therefore, geometrically, the subgroups $ G_n $ can be represented on the number line as follows

   \unitlength=0.60mm
 \begin{picture}(200,45)
 \put(0,20){\line(1,0){180}}
 \qbezier(0,20)(10,25)(20,20)
 \put(40,18){\line(0,1){4}}
 \put(60,18){\line(0,1){4}}
 \put(120,18){\line(0,1){4}}
 \qbezier(120,20)(150,10)(180,20)
 \qbezier(0,20)(30,10)(60,20)
 \put(15,9){\small ${G_{n+1}}$}
 \put(4,26){\small ${G_{n+2}}$}
 \put(55,12){\small $\frac{1}{p^{n+1}} $}
 \put(130,8){\small $ G_{n+1}\dot+(p-1)g_n$}
 \put(90,22){\small $g_n $}
 \put(90,20){\circle*{1} }
 \qbezier(0,20)(90,45)(180,20)
 \put(100,35){\small $G_n$}
 \end{picture}

\hskip4cm Figure 1. Subgroups  $G_n$.

By $X$ we denote  the collection
 of the characters of a~group $(G,\dot+ )$; this is
a~group with respect to multiplication, too. Also let
$G_n^\bot=\{\chi\in X:\forall\,x\in G_n\  , \chi(x)=1\}$ be an
annihilator of the group~$G_n$. Each annihilator~$ G_n^\bot$ is
a~group with respect to multiplication, and the subgroups~$
G_n^\bot$ form an~increa\-sing sequence
\begin{equation}
\label{eq1.03} \cdots\subset G_{-n}^\bot\subset\cdots\subset
G_0^\bot \subset G_1^\bot\subset\cdots\subset
G_n^\bot\subset\cdots
\end{equation}
with
$$
\bigcup_{n=-\infty}^{+\infty} G_n^\bot=X \quad {and} \quad
\bigcap_{n=-\infty}^{+\infty} G_n^\bot=\{1\},
$$
the quotient group $ G_{n+1}^\bot/ G_n^\bot$ having order~$p$.
The group of characters~$X$ is a zero-dimensional group with a
basic chain \ref{eq1.03}. This group  may be supplied with
topology using  a chain of subgroups~\ref{eq1.03}, the family of
the cosets $ G_n^\bot\cdot\chi$, $\chi\in X$, being taken
as~a~base of the topology. The collection of such cosets, along
with the empty set, forms the~semiring~$\mathscr{X}$. Given
a~coset $ G_n^\bot\cdot\chi$, we define a~measure~$\nu$ on it by
$\nu( G_n^\bot\cdot\chi)=\nu( G_n^\bot)= p^n$. The measure~$\nu$ can be
extended onto the $\sigma$-algebra of measurable sets in the
standard way.
 Using this measure, we form an absolutely convergent integral
$\displaystyle\int_XF(\chi)\,d\nu(\chi)$. The integral
$\displaystyle\int_G f(x)\,d\mu(x)$ is defined similarly.

The value~$\chi(g)$ of the character~$\chi$ at an element $g\in G$
will be denoted by~$(\chi,g)$. The Fourier transform~$\widehat f$
of an~$f\in L_2( G)$  is defined~as follows
$$
\widehat f(\chi)=\int_{ G}f(x)\overline{(\chi,x)}\,d\mu(x)=
\lim_{n\to+\infty}\int_{ G_{-n}}f(x)\overline{(\chi,x)}\,d\mu(x),
$$
with the limit being in the norm of $L_2(X)$. For any~$f\in L_2(G)$,
the inversion formula is valid
$$
f(x)=\int_X\widehat f(\chi)(\chi,x)\,d\nu(\chi)
=\lim_{n\to+\infty}\int_{ G_n^\bot}\widehat
f(\chi)(\chi,x)\,d\nu(\chi);
$$
here the limit also signifies the convergence in the norm of~$L_2(
G)$. If $f,g\in L_2( G)$ then the Plancherel formula is valid \cite{AVDR}
$$
\int_{ G}f(x)\overline{g(x)}\,d\mu(x)= \int_X\widehat
f(\chi)\overline{\widehat g(\chi)}\,d\nu(\chi).
$$
\goodbreak

The union of disjoint sets $E_j$ we will denote by $\bigsqcup
E_j$.
 For any $n\in \mathbb Z$ we choose a character $r_n\in  G_{n+1}^{\bot}\backslash G_n^{\bot}$
 and fixed it. The collection of functions $(r_n)_{n\in \mathbb Z}$ is called a Rademacher system. Any character $\chi$ can be rewritten as a product
 $$\chi=\prod_{j=-m}^{+\infty} r_j^{\alpha_j},\ \alpha_j=\overline{0,p-1}.$$

  Let us denote
  $$
    H_0=\{h\in G: h=a_{-1}g_{-1}\dot+a_{-2}g_{-2}\dot+\dots \dot+ a_{-s}g_{-s}, s\in \mathbb
    N,\ a_j=\overline{0,p-1}\},
  $$
  $$
    H_0^{(s)}=\{h\in G: h=a_{-1}g_{-1}\dot+a_{-2}g_{-2}\dot+\dots \dot+
    a_{-s}g_{-s},\ a_j=\overline{0,p-1}
    \},s\in \mathbb N.
  $$
  Under the Monna mapping $\lambda(H_0)= \mathbb N_0=\mathbb N\bigsqcup \{0\}$ and $\lambda(H_0^{(s)})=\mathbb N_0\bigcap [0,p^{s-1}].$ Thus the set $H_0$ is an analog of the set $\mathbb N_0$.
  \begin{definition}
 We define the mapping ${\cal A}\colon G\to G$ by
 ${\cal A}x:=\sum_{n=-\infty}^{+\infty}a_ng_{n-1}$, where
 $x=\sum_{n=-\infty}^{+\infty}a_ng_n\in G$. As any element $x\in G$ can
 be uniquely expanded as~$x=\sum a_ng_n$, the mapping ${\cal A}\colon G\to
 G$ is one-to-one onto. The mapping~${\cal A}$ is called
 a  dilation operator if~${\cal A}(x\dot+ y)={\cal A}x\dot + {\cal A}y$ for all
 $x,y\in G$.
 \end{definition}
 Since we are considering the Vilenkin group, the operator A is additive.
     By definition, put $(\chi {\cal A},x)=(\chi, {\cal
 A}x)$. It is also clear that
 ${\cal A} g_n= g_{n-1}, r_n{\cal A} = r_{n+1}$,\
 ${\cal A} G_n= G_{n-1}, G_n^\bot{\cal A}= G_{n+1}^\bot$.

  \begin{lemma}[\cite{LS2}]
 For any zero-dimensional group\\
 1) $\int\limits_{G_0^\bot}(\chi,x)\,d\nu(\chi)={\bf 1}_{G_0}(x)$,
 2) $\int\limits_{G_0}(\chi,x)\,d\mu(x)={\bf 1}_{G_0^\bot}(\chi)$.\\
   3) $\int\limits_{G_n^\bot}(\chi,x)\,d\nu(\chi)=p^n{\bf
  1}_{G_n}(x)$,
  4) $\int\limits_{G_n}(\chi,x)\,d\mu(x)=\frac{1}{p^n}{\bf
  1}_{G_n^\bot}(\chi)$.
  \end{lemma}
\begin{lemma}[\cite{LS2}]
Let  $\chi_{n,s}=r_n^{\alpha_n}r_{n+1}^{\alpha_{n+1}}\dots
r_{n+s}^{\alpha_{n+s}}$ be a character which does not belong to
$G_n^\bot$. Then
$$
\int\limits_{G_n^\bot\chi_{n,s}}(\chi,x)\,d\nu(\chi)=p^n(\chi_{n,s},x){\bf
1}_{G_n}(x).
$$
\end{lemma}
\begin{lemma}[\cite{LS2}]
Let
$h_{n,s}=a_{n-1}g_{n-1}\dot+a_{n-2}g_{n-2}\dot+\dots\dot+a_{n-s}g_{n-s}\notin
G_n$. Then
$$
\int\limits_{G_n\dot+h_{n,s}}(\chi,x)\,d\mu(x)=\frac{1}{p^n}(\chi,h_{n,s}){\bf
1}_{G_n^\bot}(\chi).
$$
\end{lemma}
 \begin{definition}[\cite{LS2}]
Let $M,N\in\mathbb N$.
We denote by  ${\mathfrak D}_{G_M}(G_{-N})$ the set of functions
 $f\in L_2(G)$ such that 1) ${\rm supp}\,f\subset G_{-N}$, and 2)
 $f$ is constant on cosets $G_M\dot+g$. The class ${\mathfrak
 D}_{G_{-N}^\bot}(G_{M}^\bot)$ is defined similarly.
 \end{definition}

 \begin{lemma}\label{lm.1.1}
 For any fixed  $\alpha_0,\alpha_1,...,\alpha_s \in \{0,1,...,p-1\}$  the set  $H_0$ is an orthonormal basis in
   $L_2(G_0^\bot r_0^{\alpha_0}r_1^{\alpha_1}...r_s^{\alpha_s})$.
 \end{lemma}
 {\it Proof.} It is known
   (\cite{LSF}) that the set  $H_0$ is  an orthonormal system in  $L_2(G_0^\bot)$.
 Moreover,  $H_0$ is an orthonormal basis in $L_2(G_0^\bot)$. Indeed, for any  $\nu\in \mathbb N$ functions
 $
 h_j=a_{-1}g_{-1}\dot+a_{-2}g_{-2}\dot+...\dot+a_{-\nu}g_{-\nu}
 $
 are constant on cosets
 \begin{equation}\label{eq1.05}
 G_{-\nu}^\bot r_{-\nu}^{\alpha_{-\nu}}r_{-\nu+1}^{\alpha_{-\nu+1}}...r_{-1}^{\alpha_{-1}}
 \end{equation}
   and orthogonal on the set    $G_0^\bot$.
    Therefore, any step function that is constant on cosets (\ref{eq1.05})
   can be uniquely represented as  linear combination
of elements
    $$
 h_j=a_{-1}g_{-1}\dot+a_{-2}g_{-2}\dot+...\dot+a_{-\nu}g_{-\nu}.
 $$
  So the system  $H_0$ is an orthonormal basis in $L_2(G_0^\bot)$.
 Let us show that $H_0$ is an orthonormal basis in  $L_2(G_0^\bot r_0^{\alpha_0}r_1^{\alpha_1}...r_\nu^{\alpha_s})$. For any $h_1,h_2\in H_0$
 $$
   \int_{G_0^\bot r_0^{\alpha_0}r_1^{\alpha_1}...r_s^{\alpha_s}}(\chi,h_1)\overline{(\chi,h_2)}d\nu(\chi)=
   \int_X{\bf 1}_{G_0^\bot r_0^{\alpha_0}r_1^{\alpha_1}...r_s^{\alpha_s}}(\chi)(\chi,h_1)\overline{(\chi,h_2)}d\nu(\chi)=
  $$
  $$
   =
   \int_X{\bf 1}_{G_0^\bot r_0^{\alpha_0}r_1^{\alpha_1}...r_s^{\alpha_s}}(\chi r_0^{\alpha_0}r_1^{\alpha_1}...r_s^{\alpha_s})(\chi r_0^{\alpha_0}r_1^{\alpha_1}...r_s^{\alpha_s},h_1)\overline{(\chi r_0^{\alpha_0}r_1^{\alpha_1}...r_s^{\alpha_s},h_2)}d\nu(\chi)=
  $$
  $$ =
   \int_X{\bf 1}_{G_0^\bot }(\chi)(\chi r_0^{\alpha_0}r_1^{\alpha_1}...r_s^{\alpha_s},h_1)\overline{(\chi r_0^{\alpha_0}r_1^{\alpha_1}...r_s^{\alpha_s},h_2)}d\nu(\chi)=
  $$
  $$ =
   (r_0^{\alpha_0}r_1^{\alpha_1}...r_s^{\alpha_s},h_1)
   \overline{(r_0^{\alpha_0}r_1^{\alpha_1}...r_s^{\alpha_s},h_2)}
   \int_X{\bf 1}_{G_0^\bot }(\chi)(\chi ,h_1)\overline{(\chi,h_2)}d\nu(\chi)=
  $$
  $$ =
   (r_0^{\alpha_0}r_1^{\alpha_1}...r_s^{\alpha_s},h_1)
   \overline{(r_0^{\alpha_0}r_1^{\alpha_1}...r_s^{\alpha_s},h_2)}
   \int_{G_0^\bot} (\chi ,h_1)\overline{(\chi,h_2)}d\nu(\chi)=\delta_{h_1,h_2}.
  $$
 This means that $H_0$ is an orthonormal system.   Let us show that $H_0$ is a basis in  $L_2(G_0^\bot r_0^{\alpha_0}r_1^{\alpha_1}...r_\nu^{\alpha_s})$.
    For any $\nu\in \mathbb N$ functions
 $
 h_j=a_{-1}g_{-1}\dot+a_{-2}g_{-2}\dot+...\dot+a_{-\nu}g_{-\nu}
 $
 are constant on cosets
  \begin{equation}\label{eq1.06}
 G_{-\nu}^\bot r_{-\nu}^{\alpha_{-\nu}}r_{-\nu+1}^{\alpha_{-\nu+1}}...r_{-1}^{\alpha_{-1}}
 r_0^{\alpha_0}r_1^{\alpha_1}...r_s^{\alpha_s}
 \end{equation}
  and orthogonal on the set   $G_0^\bot  r_0^{\alpha_0}r_1^{\alpha_1}...r_s^{\alpha_s}$.
  Therefore, any step function that is constant on cosets (\ref{eq1.06})
   can be uniquely represented as  linear combination of elements
    $$
 h_j=a_{-1}g_{-1}\dot+a_{-2}g_{-2}\dot+...\dot+a_{-\nu}g_{-\nu}.
 $$
 This means that $H_0$ is an orthonormal basis.
    $\square$

  \begin{lemma}\label{lm.1.2} Let $s\in \mathbb N$.
 For any fixed $\alpha_{-1},...,\alpha_{-s}\in \{0,1,...,p-1\}$ the family   $p^{\frac{s}{2}}\mathcal{A}^sH_0$
  is  an orthonormal basis in  $L_2(G_{-s}^\bot r_{-s}^{\alpha_{-s}}...r_{-1}^{\alpha_{-1}} )$.
 \end{lemma}
   {\it Proof.} Suppose $\tilde{h}_{1},\tilde{h}_{2}\in \mathcal{A}^s H_0$;  then
  $$
  \int_{G_{-s}^\bot r_{-s}^{\alpha_{-s}}...r_{-1}^{\alpha_{-1}}}(\chi,\tilde{h}_{1})\overline{(\chi,\tilde{h}_{2})}d\nu(\chi)=
  \int_{G_{-s}^\bot r_{-s}^{\alpha_{-s}}...r_{-1}^{\alpha_{-1}}}(\chi\mathcal{A}^s,\mathcal{A}^{-s}\tilde{h}_{1})
  \overline{(\chi\mathcal{A}^s,\mathcal{A}^{-s}\tilde{h}_{2})}d\nu(\chi)=
  $$
  $$
  \int_X{\bf 1}_{G_{-s}^\bot r_{-s}^{\alpha_{-s}}...r_{-1}^{\alpha_{-1}}}(\chi)(\chi\mathcal{A}^s,\mathcal{A}^{-s}\tilde{h}_{1})
  \overline{(\chi\mathcal{A}^s,\mathcal{A}^{-s}\tilde{h}_{2})}d\nu(\chi)=
   $$
  $$
  =\frac{1}{p^s}\int_X{\bf 1}_{G_{-s}^\bot r_{-s}^{\alpha_{-s}}...r_{-1}^{\alpha_{-1}}}(\chi\mathcal{A}^{-s})(\chi,\mathcal{A}^{-s}\tilde{h}_{1})
  \overline{(\chi,\mathcal{A}^{-s}\tilde{h}_{2})}d\nu(\chi)=
   $$
   $$
  =\frac{1}{p^s}\int_{G_{0}^\bot r_0^{\alpha_{-s}}r_1^{\alpha_{-s+1}}... r_{s-1}^{\alpha_{-1}}   }(\chi,\mathcal{A}^{-s}\tilde{h}_{1})
  \overline{(\chi,\mathcal{A}^{-s}\tilde{h}_{2})}d\nu(\chi)=\frac{1}{p^s}\delta_{\tilde{h}_1,\tilde{h}_2}.
   $$
   So, the system  $p^{\frac{s}{2}}\mathcal{A}^sH_0$ is orthonormal.   A basis property  is proved as before. $\square$

 For a given function $\varphi\in  {\mathfrak D}_{G_M}(G_{-N})$ , we define  subspaces $V_n\subset L_2(G)$ generated by $\varphi$ as
  $$
 V_n=\overline{{\rm span}\{\varphi(\mathcal{A}^n\cdot \dot-h),h\in H_0\}}, \ n\in\mathbb Z.
 $$

A function $\varphi \in L_2(G)$  is called refinable if
  \begin{equation}                                      \label{eq1.2}
   \varphi(x)=p\sum_{h\in H_0}\beta_h\varphi({\cal
   A}x\dot-h),
   \end{equation}
      for some sequence  $(\beta_h)\in l^2$. The equality   (\ref{eq1.2}) is called an refinement  equation. In Fourier domain, the equality   (\ref{eq1.2})  can be rewritten as
   \begin{equation} \label{eq1.21}
       \hat\varphi(\chi)=m_0(\chi)\hat\varphi(\chi{\cal A}^{-1}),
   \end{equation}
 where
 \begin{equation}                                           \label{eq1.3}
 m_0(\chi)=\sum_{h\in
 H_0}\beta_h\overline{(\chi{\cal A}^{-1},h)}
 \end{equation}
   is a mask of  (\ref{eq1.2}).
   \begin{lemma}[\cite{LVd2}]
    If a refinable function  $\varphi\in \mathfrak{D}_{G_{M}}(G_{-N}),\ M,N\in \mathbb N$, then its refinement equation is
      \begin{equation}                                      \label{eq1.4}
   \varphi(x)=p\sum_{h\in H_0^{(N+1)}}\beta_h\varphi({\cal
   A}x\dot-h),
   \end{equation}
   \end{lemma}

The construction of tight frame systems starts with the construction of are system $\Psi=\{\psi_1,...,\psi_q\}\subset L_2(G)$. The objective  of MRA-construction  of tight  wavelet based frames is to find $\Psi=\{\psi_1,...,\psi_q\}\subset  V_1$ such that
$$
L(\Psi):=\{\psi_{\ell,n,h}=p^\frac{n}{2}\psi_\ell(\mathcal{A}^n\cdot \dot - h):1\le\ell \le q; n\in\mathbb{Z}, h\in H_0\}
$$
 forms a tight  frame for $L_2(G)$. The system $L(\Psi)\subset L_2(G)$ is called a {\it tight wavelet frame}
of $ L_2(G)$ if
$$
\|f\|^2_{L_2(G)}=\sum\limits_{g\in L(\Psi)}|\langle f,g\rangle|^2,
$$
holds for all $f\in L_2(G)$, where $\langle \cdot,\cdot\rangle$ is the inner product in $L_2(G)$. It is equivalent to
$$
f=\sum\limits_{g\in L(\Psi)}\langle f,g\rangle g,
$$
 for all $f\in L_2(G) $.
 Since $V_1$ is a $H_0$-invariant subspace generated by $\varphi(\mathcal{A}\cdot)$, finding $\Psi\subset V_1$ is the same as finding $\beta_h^{(\ell)}$ such that
 $$
 \psi_\ell (x)=p\sum\limits_{h\in H_0}\beta_h^{(\ell)}\varphi(\mathcal{A}x\dot- h).
 $$
 In  Fourier domain, this equality     can be rewritten as
  $$
 \hat\psi_\ell(\chi)=m_\ell(\chi)\hat\varphi(\chi{\cal A}^{-1}).
 $$

We will consider only step refinable function
$\varphi(x)\in \mathfrak{D}_{G_{M}}(G_{-N})$, that is equivalent to
$\hat\varphi(\chi)\in \mathfrak{D}_{G_{-N}^\bot}(G_M^\bot),\ M,N\in\mathbb N$.

\section{Refinable functions in zero-dimensional groups}
In this article, we will use another method for constructing tight wavelet  frames other than  \cite{FLS}. First  we construct a  refinable function
 $\varphi\in \mathfrak D_{G_M}(G_{-N})$, i.e.
$\hat{\varphi}\in \mathfrak D_{G_{-N}^\bot}(G_{M}^\bot)$.
To construct such  refinable function we will use a concept of N-valid tree \cite{LSBG}. This concept  was proposed in \cite{BGLS,LBK}. A similar method for constructing scaling functions is considered in \cite{YuF6MS}.

\begin{definition}The tree T is called as N-valid if the following properties are
satisfied\\
 1) The root of this tree and its vertices of level $1,2,\dots,N-1$ are equal to zero.\\
 2) Any path $(\alpha_k\to\alpha_{k+1}\to\dots\to\alpha_{k+N-1}),\ \alpha_j=\overline{0,p-1}$
 of length $N-1$ is present in the tree $T$ exactly 1 time.
\end{definition}

 For example for $p=3, N=2$ we can construct the tree $T_4 $  of height $H=4$

\unitlength=0.75mm
\begin{picture}(100,60)

  \put(18,18){$0$}
  \put(40,20){\circle{6}}
    \put(37,20){\vector(-1,0){14}}
  \put(38,18){$0$}

   \put(52,32){\vector(-1,-1){10}}
   \put(52,8){\vector(-1,1){10}}
   \put(20,20){\circle{6}}
   \put(55,35){\circle{6}}
   \put(53,33){$2$}

   \put(55,5){\circle{6}}
   \put(53,3){$1$}

    \put(69,48){\vector(-1,-1){11}}
   \put(73,49){\circle{6}}
   \put(71,47){$0$}

   \put(70,35){\vector(-1,0){12}}
   \put(73,35){\circle{6}}
   \put(71,33){$1$}

    \put(69,22){\vector(-1,1){11}}
   \put(73,20){\circle{6}}
   \put(71,18){$2$}

   \put(99,48){\vector(-2,-1){22}}
   \put(103,50){\circle{6}}
   \put(101,48){$0$}

   \put(99,21){\vector(-2,1){22}}
   \put(103,20){\circle{6}}
   \put(101,18){$2$}

   \put(99,35){\vector(-1,0){22}}
   \put(103,35){\circle{6}}
   \put(101,33){$1$}

   \end{picture}

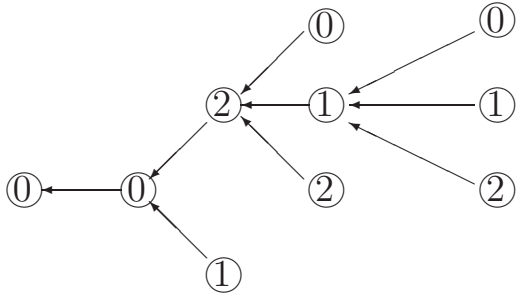
 \centerline{Figure 2. The tree $T_4$}

\noindent

Each $N$-valid tree generates a step refinable  function. Let's describe this process.
For any node $  \alpha_{-N}$ of level more than $N-1$ there exists unique path
  $$
  (  \alpha_{-1}\leftarrow\dots \leftarrow \alpha_{-N+1}\leftarrow \alpha_{-N})
  $$
  of length $N-1$ (we start numbering levels from zero). Now we  construct a graph, adding new edges to the tree $T$ in the following way.
 Let's connect  node $\alpha_{-N}$  with all pathes
 $$\alpha_{-l+N-1}\leftarrow...\leftarrow \alpha_{-l+2}\leftarrow\alpha_{-l+1}$$
  of a lower level so that
  $$\alpha_{-1}\leftarrow...\leftarrow \alpha_{-N+2}\leftarrow\alpha_{-N+1}=
  \alpha_{-l+N-1}\leftarrow...\leftarrow \alpha_{-l+2}\leftarrow\alpha_{-l+1}$$
 Denote the resulting graph by $\Gamma$.  If we add paths to the tree $T_4$, we get the graph $\Gamma_4.$

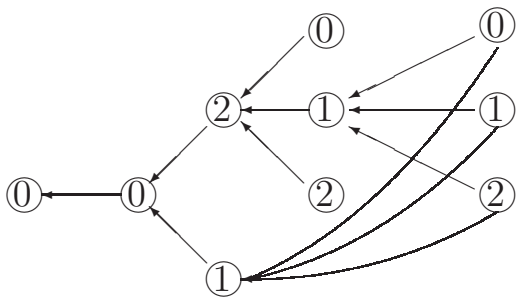
\unitlength=0.75mm
\begin{picture}(100,60)

  \put(18,18){$0$}
  \put(40,20){\circle{6}}
    \put(37,20){\vector(-1,0){14}}
  \put(38,18){$0$}

   \put(52,32){\vector(-1,-1){10}}
    \put(52,8){\vector(-1,1){10}}
   \put(20,20){\circle{6}}
   \put(55,35){\circle{6}}
   \put(53,33){$2$}

   \put(55,5){\circle{6}}
   \put(53,3){$1$}

    \put(69,48){\vector(-1,-1){11}}
   \put(73,49){\circle{6}}
   \put(71,47){$0$}

   \put(70,35){\vector(-1,0){12}}
   \put(73,35){\circle{6}}
   \put(71,33){$1$}

    \put(69,22){\vector(-1,1){11}}
   \put(73,20){\circle{6}}
   \put(71,18){$2$}

   \put(99,48){\vector(-2,-1){22}}
   \put(103,50){\circle{6}}
   \put(101,48){$0$}

   \put(99,21){\vector(-2,1){22}}
   \put(103,20){\circle{6}}
   \put(101,18){$2$}

   \put(99,35){\vector(-1,0){22}}
   \put(103,35){\circle{6}}
   \put(101,33){$1$}

\qbezier(59,5)(83,15)(103,46)
\qbezier(59,5)(83,10)(103,32)
\qbezier(59,5)(83,5)(103,17)
\put(61,5){\vector(-1,0){3}}
   \end{picture}

 \centerline{Figure 3.  Graph $\Gamma_4$ after adding new edges}

 Using vertices $
 (  \alpha_{-N},  \alpha_{-N+1},\dots , \alpha_{-1},\alpha_0)$ of the graph $\Gamma$
     we construct an $(N + 1)$-dimensional array $\lambda_{ \alpha_{-N}, \alpha_{-N+1},\dots , \alpha_{-1}, \alpha_{0}}$, enumerated by the elements of $GF(p)$ so, that \\
  (a)$\lambda_{0,0,...,0}=1$,\\
  (b)$\lambda_{ \alpha_{-N},  \alpha_{-N+1},\dots , \alpha_{-1}, \alpha_{0}}=1$ for $(\alpha_{-N},  \alpha_{-N+1},\dots , \alpha_{-1}, \alpha_{0})\in \Gamma$,\\
  (c)(b)$\lambda_{ \alpha_{-N},  \alpha_{-N+1},\dots , \alpha_{-1}, \alpha_{0}}=0$ for $(\alpha_{-N},  \alpha_{-N+1},\dots , \alpha_{-1}, \alpha_{0})\notin \Gamma$.\\
   Define mask $m_0$ on $G_1^\bot$ by the equality
   \begin{equation}                                \label{eq5.1}
  m_0(G_{-N}^\bot  r_{-N}^{\alpha_{-N}}  r_{-N+1}^{\alpha_{-N+1}}... r_{-1}^{\alpha_{-1}} r_{0}^{\alpha_{0}})=\lambda_{\alpha_{-N},  \alpha_{-N+1},\dots , \alpha_{-1}, \alpha_{0}}
   \end{equation}
   and periodically extend it on  $X$ so that
   $$
   m_0(\chi  r_{1}^{\alpha_{1}} r_{2}^{\alpha_{2}}... r_{l}^{\alpha_{l}})=m_0(\chi) $$
    for $\chi \in G_1$ and for all $ r_{1}, r_{2},..., r_{l}$.
   After that we define
   the refinable function
   \begin{equation}                                \label{eq5.2}
   \hat{\varphi}(\chi)=\prod_{k=0}^{\infty}m_0(\chi {\cal A} ^{-k}).
   \end{equation}
 In this case we say that the function $\varphi$ is generated bay the tree  $T$.
   \begin{lemma}\cite{BGLS} Let $T$ be a $N$-valid tree with height $H$. Then\\
   (a) $\hat{\varphi}(\chi)=\prod_{k=0}^{H-N+1}m_0(\chi {\cal A} ^{-k})$ for $\chi\in G_{H-2N+2}^\bot\setminus G_{H-2N+1}^\bot$,\\
   (b) $\hat{\varphi}(\chi)=0$ for $\chi\in G_{H-2N+2}^\bot\setminus G_{H-2N+1}^\bot$.\\
   (c) $\hat{\varphi}\in\mathfrak D_{G_{-N}^\bot}(G_{M}^\bot), \ M=H-2N+1$.
   \end{lemma}
Thus, every $N$-valid tree of height $H$ generates a refinable function $\varphi\in \mathfrak D_{{G}_{M}}({G}_{-N})$ with $M=H-2N+1$.

\section{Dual principle to construct tight frames   }
 In this section, we propose a method for constructing a tight  wavelet frames from a known step refinable function. For the function $\varphi\in L_2(G)$ we will use the standard notation
$$
\varphi_{n,h}=p^{\frac{n}{2}}\varphi(\mathcal{A}^n\cdot \dot-h), \quad h\in H_0, n\in \mathbb Z.
$$

Let $\varphi\in \mathfrak D_{{G}_{M}}({G}_{-N})$ be a step refinable function with a mask  $m_0({\chi})$.
  We want to find masks $m_j, j=\overline{1,q}$ and corresponding functions $\psi^{(j)}$ that generate a tight wavelet frame. We will use lemma  \ref{Lm3.1} and theorem \ref{Th3.1} which are valid for any zero-dimensional group $G$.

\begin{lemma}\label{Lm3.1}
Let the mask $m_j \ (j=1,...,q)$ satisfy the following conditions:
1)$
\hat{\psi}^{(j)}(\chi)=\hat{\varphi}(\chi\mathcal{A}^{-1})m_j(\chi)=
\xi_{\gamma_{-s},...,\gamma_{-1},\gamma_{0},...,\gamma_{u}}\in \mathbb C $   on the coset
 $$   G^\bot_{-s}r_{-s}^{\gamma_{-s}}
... r_{-1}^{\gamma_{-1}}r_{0}^{\gamma_{0}}...r_{u}^{\gamma_{u}}, \quad s=s(j), s\le N.
$$
  2) $m_j(\chi)=0$ outside the coset $   G^\bot_{-s}r_{-s}^{\gamma_{-s}}
... r_{-1}^{\gamma_{-1}}r_{0}^{\gamma_{0}}...r_{u}^{\gamma_{u}}$, \\
3)$ |\xi_{\gamma_{-s},...,\gamma_{-1},\gamma_{0},...,\gamma_{u}}|=1$.\\
Then
$$
\sum_{h\in H_0} |c_{n,h}^{(j)}(f)|^2=
\sum_{h\in H_0} |(\psi^{(j)}_{n,h},f)|^2=
\int_{G_{n-s}^\bot r_{n-s}^{\gamma_{-s}}...r_{n-1}^{\gamma_{-1}}r_{n-0}^{\gamma_{0}}...r_{n+u}^{\gamma_{u}}}|\hat{f}(\chi) |^2d\nu(\chi).
$$
\end{lemma}
{\bf Proof.} By the properties of the Fourier transform, we have
\begin{equation}\label{Eq3.1}
c_{n,h}^{(j)}(f)=
p^{\frac{n}{2}}\int_{G}f(x)\overline{\psi^{(j)}(\mathcal{A}^nx\dot-h)}d\mu(x)=
p^{\frac{n}{2}}\int_{X}\hat{f}(\chi)\overline{\hat{\psi}^{(j)}_{\mathcal{A}^n\cdot \dot-h}(\chi)}d\nu(\chi).
\end{equation}
Calculate the Fourier transform
    $$
    \hat{\psi}^{(j)}_{\mathcal{A}^n\cdot \dot-h}(\chi)=\int_{G}\psi^{(j)}(\mathcal{A}^nx \dot -h)\overline{(\chi,x)}d\mu(x)=
    $$
    $$
    \frac{1}{p^n}\int_{G}\psi^{(j)}(x \dot -h)\overline{(\chi\mathcal{A}^{-n},x)}d\mu(x)=
    $$
     $$
    =\frac{1}{p^n}\int_{G}\psi^{(j)}(x) \overline{(\chi\mathcal{A}^{-n},x\dot+ h)}d\mu(x)=
    \frac{1}{p^n}\hat{\psi}^{(j)}(\chi\mathcal{A}^{-n})\overline{(\chi\mathcal{A}^{-n},h)}.
    $$
    Substituting    in  (\ref{Eq3.1}) , we get
    $$
    c_{n,h}^{(j)}(f)=p^{-\frac{n}{2}}\int_{X}\hat{f}(\chi)
    \overline{\hat{\psi}^{(j)}(\chi\mathcal{A}^{-n})}
    (\chi\mathcal{A}^{-n},h)d\nu(\chi)=
    $$
    $$
    =p^{-\frac{n}{2}}\int_{X}\hat{f}(\chi)\overline{\hat{\varphi}(\chi\mathcal{A}^{-n-1})
    m_j(\chi\mathcal{A}^{-n})}(\chi\mathcal{A}^{-n},h)d\nu(\chi)=
    $$
    $$
    =p^{\frac{n}{2}}\int_{X}\hat{f}(\chi \mathcal{A}^{n})
    \overline{\hat{\varphi}(\chi\mathcal{A}^{-1})
    m_j(\chi)}(\chi,h)d\nu(\chi)=
    $$
    $$
    p^{\frac{n}{2}}\int_{G_{-s}^\bot r_{-s}^{\gamma_{-s}}...r_{0}^{\gamma_{0}}...r_{u}^{\gamma_{u}}  }\hat{f}(\chi \mathcal{A}^{n})\bar{\xi}_{\gamma_{-s},...,\gamma_{-1},\gamma_{0},...,\gamma_{u}}  (\chi,h)d\nu(\chi).
    $$
    Let us denote $\tilde H_0^{(s)}=\{\tilde h=a_{-s-1}g_{-s-1}\dot+...\}$,
 $h=a_{-1}g_{-1}\dot+...\dot+a_{-s}g_{-s}\dot+ \tilde h$.
  Sinse
  $(\chi, a_{-1}g_{-1}\dot+...\dot+a_{-s}g_{-s})$  is constant on any coset  $G_{-s}^\bot r_{-s}^{\gamma_{-s}}...r_{0}^{\gamma_{0}}...r_{u}^{\gamma_{u}}$ it follows by lemma  \ref{lm.1.2} that
$$
    \sum_{h\in H_0}|c_{n,h}^{(l)}(f)|^2=
$$
$$
=\sum_{a_{-1}=0}^{p-1}...\sum_{a_{-s}=0}^{p-1}\sum_{\tilde h\in \tilde H_0^{(s)}}
    |p^{\frac{n}{2}}\int_{G_{-s}^\bot r_{-s}^{\gamma_{-s}}...r_{0}^{\gamma_{0}}...r_{u}^{\gamma_{u}}  }\hat{f}(\chi \mathcal{A}^{n})(\chi,\tilde h)d\nu(\chi)|^2=
$$
$$
=p^{s+n}\sum_{\tilde h\in \tilde H_0^{(s)}}
    |\int_{G_{-s}^\bot r_{-s}^{\gamma_{-s}}...r_{0}^{\gamma_{0}}...r_{u}^{\gamma_{u}}  }\hat{f}(\chi \mathcal{A}^{n})(\chi,\tilde h)d\nu(\chi)|^2=
$$
$$
=p^n\int_{G_{-s}^\bot r_{-s}^{\gamma_{-s}}...r_{0}^{\gamma_{0}}...r_{u}^{\gamma_{u}}  }|\hat{f}(\chi \mathcal{A}^{n})|^2d\nu(\chi)=
$$
$$
=p^n\int_X {\bf 1}_{G_{-s}^\bot r_{-s}^{\gamma_{-s}}...r_{0}^{\gamma_{0}}...r_{u}^{\gamma_{u}}}(\chi)|\hat{f}(\chi \mathcal{A}^{n})|^2d\nu(\chi)=
$$
$$
=\int_X {\bf 1}_{G_{-s}^\bot r_{-s}^{\gamma_{-s}}...r_{0}^{\gamma_{0}}...r_{u}^{\gamma_{u}}}(\chi\mathcal{A}^{-n})|\hat{f}(\chi) |^2d\nu(\chi)=
$$
$$
=\int_{G_{-s}^\bot r_{-s}^{\gamma_{-s}}...r_{0}^{\gamma_{0}}..r_{u}^{\gamma_{u}}\mathcal{A}^{n}}|\hat{f}(\chi) |^2d\nu(\chi).
$$
The lemma is proved.  $\square$

\begin{theorem}\label{Th3.1}
  Let $\varphi \in \mathfrak{D}_{G_{M}}(G_{-N})$ be a refinable function with a mask $m_0$. Define masks  $m_j: j=1,2,...,q$ so that  \\
1)$\hat\varphi (\chi\mathcal{A}^{-1})m_j(\chi)=
{\bf 1}_{E_j}(\chi)$, where  $E_j=G^\bot_{-s(j)}r_{-s(j)}^{\alpha_{-s(j)}}r_{-s(j)+1}^{\alpha_{-s(j)+1}}...r_{0}^{\alpha_{0}}...r_{M}^{\alpha_{M}}$ are disjoint cosets and  $E_j \mathcal{A}^t$ are disjoint also ,\\
2)there  are integers   $t(j)\ge 0$, such  that
$$
\bigsqcup_jE_j\mathcal{A}^{t(j)}=G_{M+1}^\bot\setminus G_M^\bot.
$$
 Then functions   $
\psi^{(1)},\psi^{(2)},...,\psi^{(q)}$ generate tight wavelet frame .
\end{theorem}

{\bf Proof.} We need to prove the equality
$$
\sum_{j=1}^{q}\sum_{n(j)=-\infty}^{+\infty}\sum_{h\in H_0}|\langle\hat{\psi}^{(j)}_{n(j),h},f\rangle|^2=\int_X|\hat{f}(\chi)|^2d
\nu (\chi).
$$
From the definition of $E_j$ we have
$$
E_j\mathcal{A}^{t(j)}=G^\bot_{t(j)-s(j)}r_{t(j)-s(j)}^{\alpha_{-s(j)}}
r_{t(j)-s(j)+1}^{\alpha_{-s(j)+1}}...r_{t(j)}^{\alpha_{0}}.
$$
By the hypothesis of the theorem, we have $\bigsqcup_{j=1}^qE_j\mathcal{A}^{t(j)}= G_1^\bot\setminus G_0$.
It follows that

\begin{equation}\label{Eq4.2}
\bigsqcup_{j=1}^q \bigsqcup_{n(j)=-\infty}^{+\infty} E_j\mathcal{A}^{n(j)}= X.
\end{equation}
By lemma \ref{Lm3.1}
$$
\sum_{h\in H_0}|\langle\hat{\psi}^{(j)}_{t(j),h},f\rangle|^2=
\int_{E_j\mathcal{A}^{t(j)}}|\hat{f(\chi)}|^2d\nu(\chi)
$$
so that
$$
\sum_{j=1}^q\sum_{h\in H_0}|\langle\hat{\psi}^{(j)}_{t(j),h},f\rangle|^2=
\int_{G_1^\bot\setminus G_0^\bot}|\hat{f}(\chi)|^2d\nu(\chi).
$$
Using \ref{Eq4.2} we get
$$
\sum_{j=1}^q  \sum_{n(j)=-\infty}^{+\infty}\sum_{h\in H_0}|\langle\hat{\psi}^{(j)}_{n(j),h},f\rangle|^2=
\int_X|\hat{f}(\chi)|^2d\nu(\chi). \square
$$

To use Theorem 3.1, we must be sure that the masks $m_j$ exist. Such masks can be constructed using $N$-valid trees.
\begin{theorem}\label{Th3.1}
Let $T$ be a $N$-valid tree in which all numbers $1,2,...,p-1$ form $N$-th level.
 Let the  mask $m_0$ and the refinable  function $\varphi$ be constructed  according to this tree. Then there are functions $\psi^{(1)},\psi^{(2)},...,\psi^{(q)}$ that generate tight wavelet frame .
\end{theorem}
{\bf Proof.} In this case $\hat{\varphi}(G_{-N+2}^\bot\mathcal{A}^{-1})\neq 0$. Therefore we can take
$$
m_{j,j_1}(\chi)={\bf 1}_{E_{j,j_1}}(\chi);
{E_{j,j_1}}={G_{-N+1}^\bot r_{-N}^j r_{-N+1}^{j_1}}, j=0, j_1\in J_1; J_1\sqcup J_2=\overline{1,p-1}
$$
and
$$
m_{j,j_2}(\chi)={\bf 1}_{E_{j,j_2}}(\chi);
{E_{j,j_2}}={G_{-N}^\bot r_{-N}^j r_{-N+1}^{j_2}}, j=\overline{1,p-1}, j_2\in J_2; .
$$

\section{Approximation order }
In this section we study approximation properties of constructing frames on Vilenkin group $G$. Since $\{\psi_{n,h}^{(j)}\}$ is a tight frame then
$$\lim\limits_{\tilde{N}\to +\infty}\|f-  \sum_{n=-\infty}^{\tilde{N}}\sum_{j=1}^q\sum_{h\in H_0}\langle f,\psi_{n,h}^{(j)}\rangle \psi_{n,h}^{(j)}\|_2=0.
$$
We will find the order of this approximation.
\begin{theorem}Let the functions $\psi_1,\psi_2,...,\psi_q$ and numbers $t(j)$ be constructed as in the theorem 3.1. Let us denote $l=\max t(j)$. Then for $\tilde{N}>N$  the inequalities
$$\|f-  \sum_{n=-\infty}^{\tilde{N}}\sum_{j=1}^q\sum_{h\in H_0}\langle f,\psi_{n,h}^{(j)}\rangle \psi_{n,h}^{(j)}\|_2
\le
 $$
 $$\le p^{\frac{M}{2}}(N+1)\sum_{n=\tilde{N}+1}^\infty \left(\int_{G_{n-l+M+1}^\bot\setminus G^\bot_{n-l+M}}|\hat{f}(\chi)|^2d\nu(\chi)\right)^\frac12 $$
are satisfied.
\end{theorem}
{\bf Proof.} By the property of tight frames  we have the inequality \cite{FLS}
$$
R_{\tilde{N}}\stackrel{\text{def}}=\|f-  \sum_{n=-\infty}^{\tilde{N}} \sum_{j=1}^q\sum_{h\in H_0}\langle f,\psi_{n,h}^{(j)}\rangle \psi_{n,h}^{(j)}\|_2=
\|\sum_{n=\tilde{N}+1}^\infty \sum_{j=1}^q \sum_{h\in H_0}\langle f,\psi_{n,h}^{(j)}\rangle \psi_{n,h}^{(j)}\|_2\le .
$$
$$
 \le\sum_{n=\tilde{N}+1}^\infty \|\sum_{j=1}^q\sum_{h\in H_0}\langle f,\psi_{n,h}^{(j)}\rangle \psi_{n,h}^{(j)}\|_2.
$$
Let us denote $c^{(j)}_{n,h}:=\langle f,\psi_{n,h}^{(j)}\rangle$. Using the equality $\hat{\psi}^{(j)}_{n,h}=p^{-\frac{n}{2}}\hat{\psi}^{(j)}(\chi\mathcal{A}^{-n})\overline{(\chi \mathcal{A}^{-n},h)}$ we have
$$
\|\sum_{j=1}^q\sum_{h\in H_0}\langle f,\psi_{n,h}^{(j)}\rangle \psi_{n,h}^{(j)}\|_2^2=\|\sum_{j=1}^q\sum_{h\in H_0} c^{(j)}_{n,h} \hat{\psi}_{n,h}^{(j)}\|_2^2=
$$
$$
=\int_X|\sum_{j=1}^q\sum_{h\in H_0} c^{(j)}_{n,h} \hat{\psi}_{n,h}^{(j)}(\chi)|^2d\nu(\chi)=
$$

$$=
\int_X|\sum_{j=1}^q\sum_{h\in H_0} c^{(j)}_{n,h} p^{-\frac{n}{2}}\hat{\psi}^{(j)}(\chi\mathcal{A}^{-n})
(\overline{\chi\mathcal{A}^{-n},h)}|^2d\nu(\chi)=
$$

$$
=\frac{1}{p^n}\int_X|\sum_{j=1}^q\sum_{h\in H_0} c^{(j)}_{n,h} \hat{\psi}^{(j)}(\chi\mathcal{A}^{-n})
(\overline{\chi\mathcal{A}^{-n},h)}|^2d\nu(\chi)=
$$

$$
=\int_X|\sum_{j=1}^q\sum_{h\in H_0} c^{(j)}_{n,h} \hat{\psi}^{(j)}(\chi)
(\overline{\chi,h)}|^2d\nu(\chi)=
$$
( we used the fact that $\hat{\psi}^{(j)}={\bf 1}_{E_j}$ and  sets  $E_j$ are disjoint)
$$
=\int_X\biggl|\sum_{j=1}^q\hat{\psi}^{(j)}(\chi)\biggr|^2 \biggl|\sum_{h\in H_0} c^{(j)}_{n,h}
(\overline{\chi,h)}\biggr|^2d\nu(\chi)
=\sum_{j=1}\int_{E_j} \biggl|\sum_{h\in H_0} c^{(j)}_{n,h}
(\overline{\chi,h)}\biggr|^2d\nu(\chi).
$$
If ${E_j}\subset G_0^\bot$, then

$$
\int_{E_j} \biggl|\sum_{h\in H_0} c^{(j)}_{n,h}
(\overline{\chi,h)}\biggr|^2d\nu(\chi)\le \int_{G_0^\bot} \biggl|\sum_{h\in H_0} c^{(j)}_{n,h}
(\overline{\chi,h)}\biggr|^2d\nu(\chi)=\sum_{h\in H_0} |c^{(j)}_{n,h}|^2.
$$
If ${E_j}\subset G_1^\bot\setminus G_0^\bot$, then
$$
\int_{E_j} \biggl|\sum_{h\in H_0} c^{(j)}_{n,h}
(\overline{\chi,h)}\biggr|^2d\nu(\chi)\le \int_{G_m^\bot r_m^\alpha  } \biggl|\sum_{h\in H_0} c^{(j)}_{n,h}
(\overline{\chi,h)}\biggr|^2d\nu(\chi)=\sum_{h\in H_0} |c^{(j)}_{n,h}|^2.
$$
If ${E_j}\subset G_{m+1}^\bot\setminus G_{m}^\bot$, then
$$
\int_{E_j} \biggl|\sum_{h\in H_0} c^{(j)}_{n,h}
(\overline{\chi,h)}\biggr|^2d\nu(\chi)\le \int_{G_m^\bot r_m^\alpha  } \biggl|\sum_{h\in H_0} c^{(j)}_{n,h}
(\overline{\chi,h)}\biggr|^2d\nu(\chi)\le
$$
$$
\le\sum_{\alpha_0,\alpha_1,...,\alpha_{m-1}}\int_{G_0^\bot r_0^{\alpha_0}...r_{m-1}^{\alpha_{m-1}} r_m^\alpha  } \biggl|\sum_{h\in H_0} c^{(j)}_{n,h}
(\overline{\chi,h)}\biggr|^2d\nu(\chi)=\sum_{\alpha_0,\alpha_1,...,\alpha_{m-1}}
\sum_{h\in H_0} |c^{(j)}_{n,h}|^2=
$$
$$
=p^m\sum_{h\in H_0} |c^{(j)}_{n,h}|^2\le p^M\sum_{h\in H_0} |c^{(j)}_{n,h}|^2
$$
So
$$
\|\sum_{j=1}^q\sum_{h\in H_0}\langle f,\psi_{n,h}^{(j)}\rangle \psi_{n,h}^{(j)}\|_2^2 \le p^{M}
\sum_{j=1}^q\sum_{h\in H_0} |c^{(j)}_{n,h}|^2.
$$
By lemma \ref{Lm3.1}
$$
\sum_{h\in H_0} |c^{(j)}_{n,h}|^2=\sum_{h\in H_0}|\langle\hat{\psi}^{(j)}_{n,h},f\rangle|^2=
\int_{E_j\mathcal{A}^n}|\hat{f}(\chi)|^2d\nu(\chi).
$$
Therefore
\begin{equation}\label{Eq4.3}
\|\sum_{j=1}^q\sum_{h\in H_0}\langle f,\psi_{n,h}^{(j)}\rangle \psi_{n,h}^{(j)}\|_2^2 \le
\sum_{j=1}^q\int_{E_j\mathcal{A}^n}|\hat{f}(\chi)|^2d\nu(\chi).
\end{equation}
Let   $l=t(j)+\alpha(j)$. It's obvious that $t(j)\le N,\alpha(j)\le N$.

If $n=l$ then

$$
\sum_{j=1}^q\int_{E_j\mathcal{A}^l}|\hat{f}(\chi)|^2d\nu=
\sum_{j=1}^q\int_{E_j\mathcal{A}^{t(j)}\mathcal{A}^{\alpha(j)}}|\hat{f}(\chi)|^2d\nu\le
$$
$$
\le\sum_{j=1}^q\left(\int_{E_j\mathcal{A}^{t(j)}}|\hat{f}(\chi)|^2d\nu+
\int_{E_j\mathcal{A}^{t(j)}\mathcal{A}}|\hat{f}(\chi)|^2d\nu+...
+\int_{E_j\mathcal{A}^{t(j)}\mathcal{A}^N}|\hat{f}(\chi)|^2d\nu
\right)=
$$
$$
=\int_{G_{M+1}^\bot\setminus G_M^\bot}|\hat{f}(\chi)|^2d\nu+\int_{G_{M+2}^\bot\setminus G_{M+1}^\bot}|\hat{f}(\chi)|^2d\nu+...+\int_{G_{M+N+1}^\bot\setminus G_{M+N}^\bot}|\hat{f}(\chi)|^2d\nu.
$$
If $n>l$ then we have
$$
\sum_{j=1}^q\int\limits_{E_j\mathcal{A}^n}|\hat{f}(\chi)|^2d\nu\le \int_{G_{n-l+1}^\bot\setminus G_{n-l}^\bot}|\hat{f}(\chi)|^2d\nu+...+
\int\limits_{G_{n-l+N+M+1}^\bot\setminus G_{n-l+N+M}^\bot}|\hat{f}(\chi)|^2d\nu.
$$
Now we can write  inequality  (\ref{Eq4.3}) as
$$
\|\sum_{j=1}^q\sum_{h\in H_0}\langle f,\psi_{n,h}^{(j)}\rangle \psi_{n,h}^{(j)}\|_2\le p^\frac{M}{2}
\sum_{k=0}^N \left(\int_{G_{n-l+1+k+M}^\bot\setminus G_{n-l+k+M}^\bot}|\hat{f}(\chi)|^2d\nu(\chi)\right)^\frac12.
$$
Summing up these inequalities, we obtain for $\tilde{N}>N$
$$
\sum_{n=\tilde{N}+1}^\infty\|\sum_{j=1}^q\sum_{h\in H_0}\langle f,\psi_{n,h}^{(j)}\rangle \psi_{n,h}^{(j)}\|_2\le p^\frac{M}{2}
\sum_{n=\tilde{N}+1}^\infty\sum_{k=0}^N \left(\int_{G_{n-l+1+k+M}^\bot\setminus G_{n-l+k+M}^\bot}|\hat{f}(\chi)|^2d\nu(\chi)\right)^\frac12.
$$
Each term on the right side of this inequality is present at most $N+1$ times. Therefore
\begin{equation}\label{Eq4.4}
\sum_{n=\tilde{N}+1}^\infty \|\sum_{j=1}^q\sum_{h\in H_0}\langle f,\psi_{n,h}^{(j)}\rangle \psi_{n,h}^{(j)}\|_2\le p^{\frac{M}{2}}
(N+1)\sum_{n=\tilde{N}+1}^\infty \left(\int_{G_{n-l+M+1}^\bot\setminus G_{n-l+M}^\bot}|\hat{f}(\chi)|^2d\nu\right)^\frac12.
\end{equation}
This completes the proof.\\
{\bf Remark.} From Theorem 4.1 we can obtain an approximation estimate for functions from  Sobolev spaces.  Let's choose an increasing sequence $\gamma_n\uparrow +\infty$ so that $\sum_{n=0}^{+\infty} \frac{1}{\gamma_n}<\infty$, and $ \gamma_n=1$ for   $n\in-\mathbb{N}$. Define the function $\gamma(\chi)=\gamma_n=\gamma(\|\chi\|)$ for $\chi\in G_{n}^\bot\setminus G_{n-1}^\bot$.
Let's transform the right side in the inequality (\ref{Eq4.4}).
$$
\sum_{n=\tilde{N}+1}^\infty \left(\int_{G_{n-l+M+1}^\bot\setminus G_{n-l+M}^\bot}|\hat{f}(\chi)|^2d\nu\right)^\frac12=
\frac{1}{p^{\frac{l}{2}}}\sum_{n=\tilde{N}+1}^\infty \left(\int_{G_{n+M+1}^\bot\setminus G_{n+M}^\bot}|\hat{f}(\chi \mathcal{A}^{-l})|^2d\nu\right)^\frac12=
$$
$$
=\frac{1}{p^{\frac{l}{2}}}\sum_{n=\tilde{N}+1}^\infty \frac{1}{\gamma_{n+1}}\left(\int_{G_{n+M+1}^\bot\setminus G_{n+M}^\bot}\gamma_{n+1}^2|\hat{f}(\chi \mathcal{A}^{-l})|^2d\nu\right)^\frac12=
$$
$$=\frac{1}{p^{\frac{l}{2}}}\sum_{n=\tilde{N}+1}^\infty \frac{1}{\gamma_{n+1}}\left(\int_{G_{n+M+1}^\bot\setminus G_{n+M}^\bot}\gamma^2(\chi)|\hat{f}(\chi \mathcal{A}^{-l})|^2d\nu\right)^\frac12\le
$$
$$
\le \frac{1}{p^{\frac{l}{2}}}\sum_{n=\tilde{N}+1}^\infty\frac{1}{\gamma_{n+1}} \left(\int_X \gamma^2(\chi)|\hat{f}(\chi \mathcal{A}^{-l})|^2d\nu\right)^\frac12.
$$
Therefore
$$
R_{\tilde{N}}\le \frac{(N+1)p^{\frac{M}{2}}}{p^{\frac{l}{2}}}\sum_{n=\tilde{N}+1}^\infty\frac{1}{\gamma_{n+1}} \left(\int_X \gamma^2(\chi )|\hat{f}(\chi \mathcal{A}^{-l})|^2d\nu\right)^\frac12=
$$
$$
=(N+1)p^{\frac{M}{2}} \left(\int_X \gamma^2(\chi \mathcal{A}^{l})|\hat{f}(\chi )|^2d\nu\right)^\frac12 \sum_{n=\tilde{N}+1}^\infty\frac{1}{\gamma_{n+1}}.
$$

\noindent
 If $G$ is  Vilenkin group and $\gamma_k=p^{km}\ (k\ge 0)$ then we obtain analog of Theorem 14 from \cite{FLS}
 \begin{equation}\label{Eq4.5}
R_{\tilde{N}}\le \frac{(N+1)}{(p^m-1)p^{m{\tilde{N}+2}}} \left(\int_X (1+\|\chi\|^{m+l})^2|\hat{f}(\chi )|^2d\nu\right)^\frac12.
\end{equation}
Indeed, if $\chi\in G_k^\bot\setminus G_{k-1}^\bot$ then
$$
\gamma(\chi)=\gamma_k=
\left\{\begin{array}{lr}
p^{km}=\|\chi\|^m,& if \ k\ge 0,\\
1,&if \  k<0.\\
\end{array}\right.
$$
Therefore $\gamma(\chi)=\max(1,\|\chi\|^m)\le (1+\|\chi\|^m)$ and
$$
\int_X\gamma^2(\chi)|\hat{f}(\chi\mathcal{A}^{-l})|^2d\nu\le \int_X(1+\|\chi\|^m)^2|\hat{f}(\chi\mathcal{A}^{-l})|^2d\nu=
$$
$$
=p^l\int_X(1+\|\chi\|^{m+l})^2|\hat{f}(\chi)|^2d\nu.
$$
If we take $\gamma_k=(k+1)^{1+\varepsilon/2},\ (\varepsilon>0, k\ge 0) $ then
$$
\gamma^2(\chi)\le (1+\log_p^+\|\chi\|)^{2+\varepsilon},
$$
 and we get the inequality
  $$
R_{\tilde{N}}\le \frac{2(N+1)}{\varepsilon (1+\tilde{N})^{\varepsilon/2}} \left(\int_X \bigl(1+l+\log_p^+\|\chi\|\bigr)^{2+\varepsilon}|\hat{f}(\chi )|^2d\nu\right)^\frac12.
$$
In these inequalities $\|\chi \|= p^n$ for $\chi\in G_n^\bot\setminus G_{n-1}^\bot$, and
$$
\log_p^+\|\chi\|=
\left\{\begin{array}{lr}
\log_p\|\chi\|,& if\ \|\chi\|>1,\\
1,&if \ \|\chi\| \le 1.\\
\end{array}\right.
$$

 {\bf Acknowledgments.}This work was supported by the Russian Science Foundation
 No  22-21-00037,   https://rscf.ru/project/22-21-00037.
\vskip1cm
  \noindent


\begin{thebibliography}{99}

 \bibitem{HZ} Dong, B., Shen, Z., Zhao, H. (Editor). Las/Park City, Mathematical series, Volume 19, {\it Mathematics in Image Processing}. AMS, USA, (2013).
 \bibitem{RSh97} Ron, A., Shen, Z.   Affine system system in $L_2(\mathbb R^d)$:The analysis of the analysis operators. {\it  Journal of Functional  Analysis} 148(2): (1997), 408--447.
    \bibitem{DHRS} I.Daubechies, B.Han, A.Ron, Z.Shen. Framelets: MRA-based constructions of wavelet frames, Applied  and Computational Harmonic analysis. 14:1, (2003), 1-46.
\bibitem{ChO}   Christensen, O. {\it An Inroduction to Frames and Riesz Bases}, Boston, Birkhauser. (2003),  pp.440.
\bibitem{YuF3} Farkov, Y.  Examples of frames on the Cantor dyadic group. {\it Journal of Mathematical Sciences}. 187(1), (2012):22-34.
\bibitem{FLS}Farkov, Y., E.Lebedeva, E., Skopina, M.  Wavelet frames on Vilenkin groups and their approximation properties. {\it International Journal of Wavelets, Multiresolution and Information Processing}. 13(5), (2015):19 pages.
\bibitem{YuF4} Farkov, Yu.A.  Wavelet frames related to Walsh functions. {\it European Journal of Mathematics}. 5(1), (2019):250-267.
\bibitem{YuF6} Farkov, Yu.A.  Wavelet tight frames in Walsh analysis. {\it Annales Univ. Sci. Budapest. Sect. Comp.} 49, (2019):161-177.
\bibitem{YuF6MS} Yu. Farkov,  M. Skopina. Step wavelets on Vilenkin groups. Journal of Mathematical Sciences (2022), https://doi.org/10.1007/s10958-022-06038-w
\bibitem{LS2}  Lukomskii, S.F.  Step refinable functions and orthogonal MRA on $p$-adic Vilenkin groups. {\it JFAA},  20(1), (2014):42-65.
    \bibitem{LSF} Lukomskii, S.F. Multiresolution analysis on zero-dimensional Abelian groups and wavelets bases. {\it Sbornik: Mathematics}. 201(5), (2010):669–691.
\bibitem{AVDR}  Agaev, G.N.,   Vilenkin, N.Ya.,  Dzhafarli, G.M.,  Rubinstein, A.I. {\it Multiplicative System of Functions and Harmonic Analysis on Zero-Dimensional Groups}. Baku, (1981), pp. 180, (in Russian).
\bibitem{Mon} Monna, A.,  {\it Analyse non-Archimedienne}, New York,Springer-Verlag. (1970).
\bibitem{LSF3}  Lukomskii, S.F.  Multiresolution analysis on product of zero-dimensional Abelian groups.{\it  J. Math. Anal. Appl.} 385(2), ,(2012): 1162–1178.
\bibitem{LVd2} Lukomskii, S.F., Vodolazov, A.M. Non-Haar MRA on local Fields of positive characteristic.{\it J. Math. Anal. Appl.} 433(2), (2016): 1415-1440.
\bibitem{LSBG} S. F. Lukomskii and G. S. Berdnikov. N-Valid trees in wavelet theory on Vilenkin groups. International Journal  of Wavelets, Multiresolution  and Information Processing, Vol. 13, No. 5 (2015),
    \bibitem{LBK}  G. S. Berdnikov, Iu. S. Kruss, . S. F. Lukomskii. On orthogonal systems of shifts of scaling function on local fields of positive characteristic. Turkish Journal of Mathematics,  V.41, N.2,(2017): 244 – 253.
  \bibitem{BGLS}      Berdnikov G. S, Lukomskii S F.. Discrete orthogonal and Riesz refinable functions on local fields of positive characteristic  G. S. Berdnikov, S F. Lukomskii // European Journal of Mathe-matics. V. 6, N 4, (2020): 1505-1522.
\end{thebibliography}
 \end{document}